\documentclass[a4paper,12pt]{article}
\usepackage{amsmath}
\usepackage{latexsym}
\usepackage{url}
\usepackage{color}
\numberwithin{equation}{section}

\def\R{{\bf R}}

\def\N{{\bf N}}

\def\d{\displaystyle}
\def\e{{\varepsilon}}

\newtheorem{thm}{Theorem}[section]

\newtheorem{prop}{Proposition}[section]
\newtheorem{rem}{Remark}[section]

\title{Note on the existence of classical solutions of derivative semilinear models
for one dimensional wave equation}
\author{
Yuki Haruyama
\footnote{
Master course, Mathematical Institute,
Tohoku University,
Aoba, Sendai 980-8578, Japan.
email: yuki.haruyama.t6@dc.tohoku.ac.jp}
, Takiko Sasaki
\footnote{
Department of Mathematical Engineering, Faculty of Engineering, Musashino University,
3-3-3 Ariake, Koto-ku, Tokyo 135-8181, Japan./
Mathematical Institute, Tohoku University,
Aoba, Sendai 980-8578, Japan.
e-mail: t-sasaki@musashino-u.ac.jp.
}\quad and
Hiroyuki Takamura
\footnote{Mathematical Institute,
Tohoku University,
Aoba, Sendai 980-8578, Japan.
e-mail: hiroyuki.takamura.a1@tohoku.ac.jp.}
}
\date{{\small\it Dedicated to Professor Takayoshi Ogawa on his sixtieth birthday}
\[
\begin{array}{ll}
\mbox{\footnotesize{\bf Keywords:}}
& \mbox{\footnotesize semilinear wave equation, derivative nonlinear term,}\\
& \mbox{\footnotesize classical solution, one dimension,  lifespan}\\
\mbox{\footnotesize{\bf MSC2020:}}
& \mbox{\footnotesize primary 35L71, secondary 35B44}\\
\end{array}
\]
}
\begin{document}
\maketitle
\begin{abstract}
This note is a supplement with a new result to the review paper
by Takamura \cite{Takamura} on nonlinear wave equations in one space dimension.
We are focusing here to the long-time existence of classical solutions
of semilinear wave equations in one space dimension,
especially with derivative nonlinear terms of product-type.
Our result is an extension of the single component case,
but it is meaningful to provide models as possible as many
to cover the optimality of the general theory.
The proof is based on the classical iteration argument of the point-wise estimate of the solution.
\end{abstract}


\section{Introduction}

\par
In this note, we shall consider the initial value problems;
\begin{equation}
\label{IVP}
\left\{
\begin{array}{ll}
	\d u_{tt}-u_{xx}=|u_t|^p|u_x|^q
	&\mbox{in}\quad \R\times(0,T),\\
	u(x,0)=\e f(x),\ u_t(x,0)=\e g(x),
	& x\in\R,
\end{array}
\right.
\end{equation}
where $p,q\in (1,\infty)\cup\{0\}$ and $T>0$.
Let us assume that $f$ and $g$ are given smooth functions of compact support
and a parameter $\e>0$ is \lq\lq small enough".
We are interested in the lifespan $T(\e)$, the maximal existence time,
of classical solutions of (\ref{IVP}).
Note that we do not investigate the case of $p\in(0,1)$ or $q\in(0,1)$
because we cannot expect the existence of local-in-time solution.
It is naturally conjectured that the estimate,
\begin{equation}
\label{lifespan}
T(\e)\sim C\e^{-(p+q-1)},
\end{equation}
holds because 
the valance of the amplitude of the solution in both sides in (\ref{IVP}) is
$\e\approx\e^{p+q}t$ according to the integral representations of $u_t$ and $u_x$ in (\ref{u_t}) and (\ref{u_x})
as well as definitions of integral operators $L'$ and $\overline{L'}$ in (\ref{nonlinear_derivative}) and (\ref{nonlinear_derivative_conjugate}) below.
Here we denote the fact that there are positive constants,
$C_1$ and $C_2$, independent of $\e$ satisfying $A(\e,C_1)\le T(\e)\le A(\e,C_2)$
by $T(\e)\sim A(\e,C)$.

\par
The case of $p>1$ and $q=0$ in (\ref{lifespan}) is verified by Zhou \cite{Zhou01} for the upper bound
and by Kitamura, Morisawa and Takamura \cite{KMT23} for the lower bound. 
We note that this case was first studied to show the optimality of the general theory
for nonlinear wave equations by Li, Yu and Zhou \cite{LYZ91, LYZ92} 
in which only the smooth nonlinear terms are considered.
On the other hand, the case of $p=0$ and $q>1$ in (\ref{lifespan})
is verified by Sasaki, Takamatsu and Takamura \cite{STT23}.
This result is meaningful especially in the blow-up part,
not only for the optimality but also for the future study of \lq\lq the blow-up boundary"
because the case of $p>1$ and $q=0$ is well-studied by Sasaki \cite{Sasaki18}
but this case is still open due to the fact that there is no point-wise positivity of the solution. 
We note that \cite{STT23} has a trivial error in the proof of the existence part.
It is easy to modify it, but here we show how to do among the proof of our result.
Finally, we point out that Haruyama and Takamura \cite{HT} recently study a special quasilinear model,
$u_{tt}-u_{xx}=|u_x|^{p-2}u_xu_{xx}$ with $p>1$,
and prove that the upper bound in (\ref{lifespan}) with $p>1$ and $q=0$ is valid even for this equation.
See Introduction in \cite{HT} and references therein.

\par
Our purpose is to show the lower bound of (\ref{lifespan}) in the general case of $p>1$ and $q>1$.
We note that the upper bound in such a case is still open because the proof of the blow-up result
for the case of $p>1$ and $q=0$ is organized
by comparison argument with ordinary differential equations
and point-wise estimates of the solution itself,
while the one for the case of $p=0$ and $q>1$ is organized by ordinary differential inequality
for some weighted functional of the solution which was introduced
by Rammaha \cite{Rammaha95, Rammaha97}.
See \cite{KMT23, STT23} for details.
The difficulty in proving the blow-up part in (\ref{lifespan}) comes from this difference.
Moreover, we investigate the upper bound of the lifespan of a solution of a special model,
\begin{equation}
\label{eq:special}
u_{tt}-u_{xx}=|u_t\pm u_x|^{p-1}(u_t\pm u_x)\ (p>1),
\end{equation}
which includes product terms of $u_t$ and $u_x$.
The result is the expected one, $T(\e)\le C\e^{-(p-1)}$.
We note that this is optimal when $p$ is odd integer, or sufficiently large,
by the general theory \cite{LYZ91, LYZ92}.
\par
This paper is organized as follows.
The preliminaries and our main result are stated in Section 2.
Section 3 is devoted to the proof of the main result
by employing the classical iteration argument of point-wise estimates of the solution
which was firstly introduced by John \cite{John79}.
Finally we investigate the upper bound of the lifespan of the solution to (\ref{eq:special}) in Section 4.


\section{Preliminaries and main result}

Throughout this paper, we assume that the initial data
$(f,g)\in C_0^2(\R)\times C^1_0(\R)$ satisfies
\begin{equation}
\label{supp_initial}
\mbox{\rm supp }f,\ \mbox{supp }g\subset\{x\in\R:|x|\le R\},\quad R\ge1.
\end{equation}
Let $u$ be a classical solution of (\ref{IVP}) in the time interval $[0,T]$.
Then, due to the finiteness of the propagation speed of the wave,
the support condition of the initial data, (\ref{supp_initial}), implies that
\begin{equation}
\label{support_sol}
\mbox{supp}\ u(x,t)\subset\{(x,t)\in\R\times[0,T]:|x|\le t+R\}.
\end{equation}
For example, see Appendix of John \cite{John_book} for this fact.

\par
It is well-known that $u$ satisfies the following integral equation.
\begin{equation}
\label{u}
u(x,t)=\e u^0(x,t)+L(|u_t|^p|u_x|^q)(x,t),
\end{equation}
where $u^0$ is a solution of the free wave equation with the same initial data,
\begin{equation}
\label{u^0}
u^0(x,t):=\frac{1}{2}\{f(x+t)+f(x-t)\}+\frac{1}{2}\int_{x-t}^{x+t}g(y)dy,
\end{equation}
and a linear integral operator $L$ for a function $U=U(x,t)$ in Duhamel's term is defined by
\begin{equation}
\label{nonlinear}
L(U)(x,t):=\frac{1}{2}\int_0^tds\int_{x-t+s}^{x+t-s}U(y,s)dy.
\end{equation}
Then, one can apply the time-derivative to (\ref{u}) to obtain
\begin{equation}
\label{u_t}
u_t(x,t)=\e u_t^0(x,t)+L'(|u_t|^p|u_x|^q)(x,t)
\end{equation}
and
\begin{equation}
\label{u^0_t}
u_t^0(x,t)=\frac{1}{2}\{f'(x+t)-f'(x-t)+g(x+t)+g(x-t)\},
\end{equation}
where $L'$ for a function $U=U(x,t)$ is defined by
\begin{equation}
\label{nonlinear_derivative}
L'(U)(x,t):=\frac{1}{2}\int_0^t\{U(x+t-s,s)+U(x-t+s,s)\}ds.
\end{equation}
On the other hand, applying the space-derivative to (\ref{u}),
we have
\begin{equation}
\label{u_x}
u_x(x,t)=\e u_x^0(x,t)+\overline{L'}(|u_t|^p|u_x|^q)(x,t)
\end{equation}
and
\begin{equation}
\label{u^0_x}
u_x^0(x,t)=\frac{1}{2}\{f'(x+t)+f'(x-t)+g(x+t)-g(x-t)\},
\end{equation}
where $\overline{L'}$ for a function $U=U(x,t)$ is defined by
\begin{equation}
\label{nonlinear_derivative_conjugate}
\overline{L'}(U)(x,t):=
\frac{1}{2}\int_0^t\{U(x+t-s,s)-U(x-t+s,s)\}ds.
\end{equation}
Therefore, both $u_t$ and $u_x$ are expressed by themselves.
Moreover, one more space-derivative to (\ref{u_t}) yields that
\begin{equation}
\label{u_tx}
\begin{array}{ll}
u_{tx}(x,t)=&\d\e u_{tx}^0(x,t)\\
&\d+L'(p|u_t|^{p-2}u_tu_{tx}|u_x|^q+q|u_x|^{q-2}u_xu_{xx}|u_t|^p)(x,t)
\end{array}
\end{equation}
and
\begin{equation}
\label{u^0_tx}
u_{tx}^0(x,t)=\frac{1}{2}\{f''(x+t)-f''(x-t)+g'(x+t)+g'(x-t)\}.
\end{equation}
Similarly, we have that
\[
\begin{array}{ll}
u_{tt}(x,t)=&
\d\e u_{tt}^0(x,t)+|u_t(x,t)|^p|u_x(x,t)|^q\\
&\d +\overline{L'}(p|u_t|^{p-2}u_tu_{tx}|u_x|^q+q|u_x|^{q-2}u_xu_{xx}|u_t|^p)(x,t)
\end{array}
\]
and
\[
u_{tt}^0(x,t)=\frac{1}{2}\{f''(x+t)+f''(x-t)+g'(x+t)-g'(x-t)\}.
\]
Therefore, $u_{tt}$ is expressed by $u_t,u_x,u_{tx},u_{xx}$.
Finally we see that
\begin{equation}
\label{u_xx}
\begin{array}{ll}
u_{xx}(x,t)=
&\d\e u_{xx}^0(x,t)\\
&\d +\overline{L'}(p|u_t|^{p-2}u_tu_{tx}|u_x|^q+q|u_x|^{q-2}u_xu_{xx}|u_t|^p)(x,t)
\end{array}
\end{equation}
and
\begin{equation}
\label{u^0_xx}
u_{xx}^0(x,t)=u^0_{tt}(x,t).
\end{equation}

\par
First, we note the following fact.

\begin{prop}
\label{prop:system}
Assume that $(f,g)\in C^2(\R)\times C^1(\R)$.
Let $(v,w)$ be a $C^1$ solution of a system of integral equations;
\begin{equation}
\label{system}
\left\{
\begin{array}{l}
v=\e u_t^0+L'(|v|^p|w|^q),\\
w=\e u_x^0+\overline{L'}(|v|^p|w|^q)
\end{array}
\right.
\mbox{in}\ \R\times[0,T]
\end{equation}
with some $T>0$.
Then, $u$ defined by
\[
u(x,t)=\e u^0(x,t)+L(|v|^p|w|^q)(x,t)
\]
is a classical solution of (\ref{IVP}) in $\R\times[0,T]$,
and $v\equiv u_t, w\equiv u_x$ in $\R\times[0,T]$ hold.
\end{prop}
\par\noindent
{\bf Proof.} It is trivial that $v\equiv u_t, w\equiv u_x$ by differentiation with respect to $t$ and $x$.
 The rest part is easy along with the computations above in this section. 
\hfill$\Box$

\vskip10pt
Our main result is the following theorem.

\begin{thm}
\label{thm:lower-bound}
Let $p>1$ and $q>1$.
Assume (\ref{supp_initial}).
Then, there exists a positive constant $\e_0=\e_0(f,g,p,q,R)>0$ such that
a classical solution $u\in C^2(\R\times[0,T])$ of (\ref{IVP}) exists
as far as $T$ satisfies
\begin{equation}
\label{lower-bound}
T\le c\e^{-(p+q-1)}
\end{equation}
where $0<\e\le\e_0$, and $c$ is a positive constant independent of $\e$.
\end{thm}

\begin{rem}
\label{rem:lifespan}
Due to the definition of the lifespan,
\[
\begin{array}{ll}
T(\e):=\sup\{T>0:&\mbox{There exists a classical solution of (\ref{IVP}) in $[0,T]$}\\
&\mbox{for arbitrarily fixed data $(f,g)$. \}},
\end{array}
\]
the result of Theorem \ref{thm:lower-bound} above gives us the lower bound of the lifespan
with a statement that there exists
a positive constant $\e_0=\e_0(f,g,p,q,R)>0$ such that
\[
c\e^{-(p+q-1)}\le T(\e),
\]
where $0<\e\le\e_0$, and $c$ is a positive constant independent of $\e$.
\end{rem}


\section{Proof of Theorem \ref{thm:lower-bound}}
\par
The proof of our main result is based on the argument in Kitamura, Morisawa and Takamura \cite{KMT23}.
According to Proposition \ref{prop:system},
we shall construct a $C^1$ solution of (\ref{system}).
Let $\{(v_j,w_j)\}_{j\in\N}$ be a sequence of $\{C^1(\R\times[0,T])\}^2$ defined by
\begin{equation}
\label{v_j,w_j}
\left\{
\begin{array}{ll}
v_{j+1}=\e u_t^0+L'(|v_j|^p|w_j|^q), & v_1=\e u_t^0,\\
w_{j+1}=\e u_x^0+\overline{L'}(|v_j|^p|w_j|^q), &  w_1=\e u_x^0.
\end{array}
\right.
\end{equation}
Then, in view of (\ref{u_tx}) and (\ref{u_xx}), $\left((v_j)_x,(w_j)_x\right)$ has to satisfy
\begin{equation}
\label{v_j_x,w_j_x}
\left\{
\begin{array}{ll}
(v_{j+1})_x&=\e u_{tx}^0+L'\left(p|v_j|^{p-2}v_j(v_j)_x|w_j|^q\right)\\
&\quad+L'\left(q|w_j|^{q-2}w_j(w_j)_x|v_j|^p\right),\\
(v_1)_x&=\e u_{tx}^0,\\
(w_{j+1})_x&=\e u^0_{xx}+\overline{L'}\left(p|v_j|^{p-2}v_j(v_j)_x|w_j|^q\right)\\
&\quad+\overline{L'}\left(q|w_j|^{q-2}w_j(w_j)_x|v_j|^p\right),\\
 (w_1)_x&=\e u_{xx}^0,
\end{array}
\right.
\end{equation}
so that the function space in which $\{(v_j,w_j)\}$ converges is
\[
\begin{array}{ll}
X:=&\{(v,w)\in\{C^1(\R\times[0,T])\}^2\ :\ \|(v,w)\|_X<\infty,\\
&\quad\mbox{supp}\ (v,w)\subset\{(x,t)\in\R\times[0,T]\ :\ |x|\le t+R\}\},
\end{array}
\]
which is equipped with a norm
\[
\|(v,w)\|_X:=\|v\|+\|v_x\|+\|w\|+\|w_x\|,
\]
where
\[
\|v\|:=\sup_{(x,t)\in\R\times[0,T]}|v(x,t)|.
\]
First we note that supp $(v_j,w_j)\subset\{(x,t)\in\R\times[0,T]\ :\ |x|\le t+R\}$ implies supp
$(v_{j+1},w_{j+1})\subset\{(x,t)\in\R\times[0,T]\ :\ |x|\le t+R\}$.
It is easy to check this fact by assumption on the initial data (\ref{supp_initial})
and the definitions of $L',\overline{L'}$ in the previous section.

\par
The following lemma is a priori estimate.
\begin{prop}
\label{prop:apriori}
Let $(v,w)\in\{C(\R\times[0,T])\}^2$ and supp\ $(v,w)\subset\{(x,t)\in\R\times[0,T]:|x|\le t+R\}$. Then there exists a positive constant $C$ independent of $T$ and $\e$ such that
\begin{equation}
\label{apriori}
\begin{array}{l}
\|L'(|w|^p|u|^q)\|\le C\|w\|^p\|u\|^q(T+R),\\
\|\overline{L'}(|w|^p|u|^q)\|\le C\|w\|^p\|u\|^q(T+R).
\end{array}
\end{equation}
\end{prop}

\par\noindent
{\bf Proof.} The proof of the first inequality in (\ref{apriori}) is exactly same as
the one of Proposition 3.1 with $a=-1$ in Morisawa, Kitamura and Takamura \cite{KMT23}.
Because the total power in (\ref{apriori}) is $(p+q)$
which is equivalent to $p$ in \cite{KMT23}.
The second inequality is readily follows from the trivial inequality,
\[
|\overline{L'}(U)|\le|L'(U)|
\]
for a positive function $U=U(x,t)$.
\hfill$\Box$

\vskip10pt
\par
Once we have a priori estimate (\ref{apriori}) in Proposition \ref{prop:apriori},
the proof of Theorem \ref{thm:lower-bound} is immediately established.
Because the boundedness of the sequence in $X$ is the same
as Morisawa, Kitamura and Takamura \cite{KMT23} and
the condition by estimate of the convergence of the sequence in $X$ is harmless
if one follows Kido, Sasaki, Takamatsu and Takamura \cite{KSTT24}
to handle the product nonlinear terms.
The proof of Theorem \ref{thm:lower-bound} is now completed.
\hfill$\Box$

\vskip10pt
\par
Finally we shall state some remarks below.

\begin{rem}
In our future work, the optimality of Theorem \ref{thm:lower-bound} will be provided hopefully.
But as stated in Introduction, we will meet a technical difficulty to end.
A combination of the functional method and the point-wise estimate may help us.
\end{rem}

\begin{rem}
As promised in Introduction, we mention to an erratum to Sasaki, Takamatsu and Takamura \cite{STT23}.
In \cite{STT23},  there is a trivial error to define the function space for Theorem  2.1.
The norm of $X$ in \cite{STT23} should be
\[
\|w\|_X=\|w\|+\|w_x\|
\]
like this paper.
With this correction, $L'$ in the proof of Theorem 2.1 in \cite{STT23}
should be replaced with $\overline{L'}$.
But it provides nothing new again by $|\overline{L'}(U)|\le|L'(U)|$ for a positive function $U=U(x,t)$.
\end{rem}


\section{Blow-up example including product terms}
\par
As promised in the introduction, we shall prove the following theorem.

\begin{thm}
\label{thm:special}
Let $f\in C_0^2(\R)$ and $g\in C_0^1(\R)$ satisfy
\begin{equation}
\label{positive}
\pm f'(x_0)+g(x_0)\neq0
\end{equation}
for some $x_0\in\R$.
If the equation in (\ref{IVP}) is replaced with the one in (\ref{eq:special}),
then there exists a positive constant $C$ such that
a classical solution cannot exists if $T$ satisfies
\begin{equation}
\label{lifespan_special}
T>C\e^{-(p-1)}
\end{equation}
for all $\e>0$.
\end{thm}

\begin{rem}
(\ref{lifespan_special}) implies the lifespan estimate
\[
T(\e)\le C\e^{-(p-1)}.
\]
When $p$ is odd integer, or sufficiently large, this is known to be sharp by general theory \cite{LYZ91, LYZ92}.
\end{rem}

\begin{rem}
If the blow-up condition (\ref{positive}) is false, then we may have a global-in-time solution
with arbitrary $\e>0$.
In fact, it is easy to see that $u(x,t):=\e f(x\mp t)$ solves the equation (\ref{eq:special})
for which the initial condition should be formed as
\[
u(x,0)=\e f(x),\ u_t(x,0)=\mp\e f'(x),
\]
so that it makes the assumption,
\[
\pm f'(x)+g(x)=0\quad\mbox{for all }x\in\R.
\]

\end{rem}

\par\noindent
{\bf Proof of Theorem \ref{thm:special}.}
Let $u$ be a classical solution of this problem.
It follows from (\ref{u_t}) and (\ref{u_x}) that
\begin{equation}
\label{u_tpmu_x}
\begin{array}{ll}
u_t(x,t)\pm u_x(x,t)
&\d=\e\{\pm f'(x\pm t)+g(x\pm t)\}\\
&\d\quad+\int_0^t|u_t\pm u_x|^{p-1}(u_t\pm u_x)(x\pm(t-s),s)ds.
\end{array}
\end{equation}
First, let us assume that
\begin{equation}
\label{data}
M_\pm:=\pm f'(x_0)+g(x_0)>0
\end{equation}
and fix
\[
x\pm t=x_0.
\]
Moreover, set
\[
U_\pm(t):=u_t(x_0\mp t,t)\pm u_x(x_0\mp t,t)
\]
Then, we can rewrite (\ref{u_tpmu_x}) as
\[
U_\pm(t)=M_\pm\e
+\int_0^t|U_\pm(s)|^{p-1}U_\pm(s)ds.
\]
Hence the continuity of $U_\pm$ and the positivity of the data $U_\pm(0)=M_\pm\e>0$ yield that
\[
U_\pm(t)>0\quad\mbox{for}\ t\ge0.
\]
Therefore $U_\pm$ satisfied
\[
\left\{
\begin{array}{l}
U_\pm'=U_\pm^p\\
U_\pm(0)=M_\pm\e
\end{array}
\right.
\]
which is solvable as
\[
U_\pm(t)=\{(M_\pm\e)^{1-p}-(p-1)t\}^{-1/(p-1)}
\]
and has an exact blow-up time
\begin{equation}
\label{blow-up_time}
t_0:=\frac{(M_\pm\e)^{1-p}}{p-1}>0.
\end{equation}

\par
When $M_{\pm}<0$ instead of (\ref{data}), then one has the same result as above by
substituting $u(x,t)$ by $-u(x,t)$.
Therefore, Theorem \ref{thm:special} is now established with $C=|M_\pm|^{1-p}/(p-1)>0$.
\hfill$\Box$

\begin{rem}
According to the analysis on the blow-up boundary of solutions of
transport equations in Takeno \cite{Takeno},
we may have a similar result also on our equation (\ref{eq:special})
in view of the concrete expression of the blow-up time $t_0$ in (\ref{blow-up_time})
which is a relation between $t_0$ and $x_0$.
\end{rem}

\section*{Acknowledgement}
\par
The second and third authors are partially supported
by the Grant-in-Aid for Scientific Research (A) (No.22H00097) and (C) (No.24K06819), 
Japan Society for the Promotion of Science.
All the authors are grateful to the anonymous referee for his/her useful suggestions
to complete the manuscript


\bibliographystyle{plain}

\end{document}